\documentclass[a4paper, 11pt]{article}
\usepackage[ansinew]{inputenc}
\usepackage{amssymb, amsmath, amsthm}
\usepackage{nicefrac, mathdots}

\def\typeout#1{\message{^^J}\message{#1}\message{^^J}}
\newif\ifSRCOK \SRCOKtrue
\newcount\PAGETOP
\newcount\LASTLINE
\global\PAGETOP=1
\global\LASTLINE=-1
\def\EJECT{\SRC\eject}
\def\WinEdt#1{\typeout{:#1}}
\gdef\MainFile{\jobname.tex}
\gdef\CurrentInput{\MainFile}
\newcount\INPSP
\global\INPSP=0
\def\SRC{\ifSRCOK%
  \ifnum\inputlineno>\LASTLINE%
    \ifnum\LASTLINE<0%
      \global\PAGETOP=\inputlineno%
    \fi%
    \global\LASTLINE=\inputlineno%
    \ifnum\INPSP=0%
      \ifnum\inputlineno>\PAGETOP%
        
      \fi%
    \else%
      
    \fi%
  \fi%
\fi}
\def\PUSH#1{%
\SRC%
\ifnum\INPSP=0 \global\let\INPSTACKA=\CurrentInput \else%
\ifnum\INPSP=1 \global\let\INPSTACKB=\CurrentInput \else%
\ifnum\INPSP=2 \global\let\INPSTACKC=\CurrentInput \else%
\ifnum\INPSP=3 \global\let\INPSTACKD=\CurrentInput \else%
\ifnum\INPSP=4 \global\let\INPSTACKE=\CurrentInput \else%
\ifnum\INPSP=5 \global\let\INPSTACKF=\CurrentInput \else%
               \global\let\INPSTACKX=\CurrentInput \fi\fi\fi\fi\fi\fi%
\gdef\CurrentInput{#1}%
\WinEdt{<+ \CurrentInput}%
\global\LASTLINE=0%
\ifSRCOK\fi%
\global\advance\INPSP by 1}
\def\POP{%
\ifnum\INPSP>0 \global\advance\INPSP by -1  \fi%
\ifnum\INPSP=0 \global\let\CurrentInput=\INPSTACKA \else%
\ifnum\INPSP=1 \global\let\CurrentInput=\INPSTACKB \else%
\ifnum\INPSP=2 \global\let\CurrentInput=\INPSTACKC \else%
\ifnum\INPSP=3 \global\let\CurrentInput=\INPSTACKD \else%
\ifnum\INPSP=4 \global\let\CurrentInput=\INPSTACKE \else%
\ifnum\INPSP=5 \global\let\CurrentInput=\INPSTACKF \else%
               \global\let\CurrentInput=\INPSTACKX \fi\fi\fi\fi\fi\fi%
\WinEdt{<-}%
\global\LASTLINE=\inputlineno%
\global\advance\LASTLINE by -1%
\SRC}
\def\INPUT#1{\relax}

\let\originalxxxeverypar\everypar
\newtoks\everypar
\originalxxxeverypar{\the\everypar\expandafter\SRC}
\everymath\expandafter{\the\everymath\expandafter\SRC}
\output\expandafter{\expandafter\SRCOKfalse\the\output}

\newif\ifSRCOK \SRCOKtrue
\newcount\PAGETOP
\newcount\LASTLINE
\global\PAGETOP=1
\global\LASTLINE=-1
\gdef\MainFile{\jobname.tex}
\gdef\CurrentInput{\MainFile}
\newcount\INPSP
\global\INPSP=0
\def\EJECT{\SRC\eject}
\def\WinEdt#1{\typeout{:#1}}
\def\SRC{\ifSRCOK%
  \ifnum\inputlineno>\LASTLINE%
    \ifnum\LASTLINE<0%
      \global\PAGETOP=\inputlineno%
    \fi%
    \global\LASTLINE=\inputlineno%
    \ifnum\INPSP=0%
      \ifnum\inputlineno>\PAGETOP%
      \fi%
    \else%
    \fi%
  \fi%
\fi}
\def\PUSH#1{%
\SRC%
\ifnum\INPSP=0 \global\let\INPSTACKA=\CurrentInput \else%
\ifnum\INPSP=1 \global\let\INPSTACKB=\CurrentInput \else%
\ifnum\INPSP=2 \global\let\INPSTACKC=\CurrentInput \else%
\ifnum\INPSP=3 \global\let\INPSTACKD=\CurrentInput \else%
\ifnum\INPSP=4 \global\let\INPSTACKE=\CurrentInput \else%
\ifnum\INPSP=5 \global\let\INPSTACKF=\CurrentInput \else%
               \global\let\INPSTACKX=\CurrentInput \fi\fi\fi\fi\fi\fi%
\gdef\CurrentInput{#1}%
\WinEdt{<+ \CurrentInput}%
\global\LASTLINE=0%
\ifSRCOK\fi%
\global\advance\INPSP by 1}
\def\POP{%
\ifnum\INPSP>0 \global\advance\INPSP by -1  \fi%
\ifnum\INPSP=0 \global\let\CurrentInput=\INPSTACKA \else%
\ifnum\INPSP=1 \global\let\CurrentInput=\INPSTACKB \else%
\ifnum\INPSP=2 \global\let\CurrentInput=\INPSTACKC \else%
\ifnum\INPSP=3 \global\let\CurrentInput=\INPSTACKD \else%
\ifnum\INPSP=4 \global\let\CurrentInput=\INPSTACKE \else%
\ifnum\INPSP=5 \global\let\CurrentInput=\INPSTACKF \else%
               \global\let\CurrentInput=\INPSTACKX \fi\fi\fi\fi\fi\fi%
\WinEdt{<-}%
\global\LASTLINE=\inputlineno%
\global\advance\LASTLINE by -1%
\SRC}
\def\INPUT#1{\relax}
\let\OldINCLUDE=\include
\def\include#1{
\EJECT%
\PUSH{#1.tex}%
\OldINCLUDE{#1}%
\POP}

\let\originalxxxeverypar\everypar
\newtoks\everypar
\originalxxxeverypar{\the\everypar\expandafter\SRC}
\everymath\expandafter{\the\everymath\expandafter\SRC}
\let\zzzxxxbibliography=\bibliography
\def\bibliography#1{\PUSH{\jobname.bbl}\zzzxxxbibliography{#1}\POP}
\output\expandafter{\expandafter\SRCOKfalse\the\output}

\title{Was ist Unendlichkeit - und wenn ja, wie viele?}
\author{Martin Meyries \vspace{0.2cm}\\Institut f\"ur Mathematik\\Martin-Luther-Universit\"at Halle-Wittenberg\vspace{0.2cm}\\\texttt{martin.meyries@mathematik.uni-halle.de}}

\begin{document}

\maketitle


\subsection*{Vorbemerkung}
Georg Cantor (1845 -- 1918, ab 1872 Professor in Halle) hat den Mathe\-matikern fast im Alleingang das Tor zur Unendlichkeit ge\"offnet und durch seine daraus resultierende Mengenlehre die Mathematik von Grund auf ver\-\"an\-dert. Es sind haupts\"achlich seine Ideen, die hier pr\"asentiert werden. 

Der vorliegende Text gibt eine elementare Einf\"uhrung in die Eigenschaften unendlicher Mengen. Der Zugang soll dabei so einfach wie m\"oglich gemacht werden. In diesem Sinne ben\"otigt man -- im Prinzip -- keine mathe\-matischen Vorkenntnisse, unter anderem wird fast vollst\"andig auf Formeln verzichtet. Was man aber ganz sicher ben\"otigt, ist gen\"ugend Ausdauer und den Willen, auch m\"oglicherweise schwierige Dinge zu verstehen. Dies ist kein wissenschaftlicher Text, er versucht aber die wichtigsten Eigenschaften eines mathematischen Aufsatzes zu erf\"ullen: Transparenz und Rigorosit\"at. 

Behandelt werden unter anderem potentielle und aktuale Unendlichkeit, die Begriffe M\"achtigkeit und Abz\"ahlbarkeit, Hilberts Hotel, Cantors Diagonalargumente und die mehr als erstaunliche Antwort auf die Frage nach der G\"ultigkeit von Cantors Kontinuumshypothese.\medskip

\begin{small} Ich hoffe, den genannten Anspr\"uchen einigerma\ss{}en gerecht geworden zu sein und freue mich \"uber Lob, Kritik und Hinweise sowie Vorschl\"age aller Art, wie man die Thematik f\"ur Nichtmathematiker noch leichter zug\"anglich machen k\"onnte.

Ganz herzlich bedanken m\"ochte ich mich bei Alessa Binder, Andreas Bolleyer und Philip M\"uller f\"ur die kritische Durchsicht des Textes.\end{small}

\subsection*{Potentielle Unendlichkeit}
Manche Kinder behaupten, dass es eine gr\"o\ss{}te Zahl g\"abe. Wenn man sie dann fragt, welche das denn sei, erh\"alt man Antworten in den Gr\"o\ss{}en\-ord\-nungen von Hunderten, Tausenden oder Millionen, je nach Kenntnis\-stand des Kindes. Es ist nun aber eine so einfache wie fundamentale Erkenntnis, dass es keine gr\"o\ss{}te Zahl gibt, da man ja zum Beispiel auf jede Zahl immer noch 1 dazuaddieren kann.

 Im Alltag kommt man mit Zahlen bis 10.000 gut aus und kann damit meist auch noch konkrete Dinge verbinden (vielleicht ist man 1000 Kilometer weit mit dem Auto in den Urlaub gefahren, hat 950 Gramm Mehl abgewogen oder 10.000 Euro BAf\"oG zur\"uckgezahlt). Schon viel schwerer vorstellbar sind die Millionen, Milliarden oder Billionen, mit denen in Politik, Wirtschaft und Medien jongliert wird. 
 
Noch gr\"o\ss{}ere Zahlen, die eine reale Bedeutung haben, findet man meist nur in Naturwissenschaften und Technik. Beispielsweise wiegt die Erde ca. 
$6.000.000.000.000.000.000.000.000$  (sechs Quadrillionen) Kilogramm, und die Sonne ist sogar noch eine Million mal schwerer. Unvorstellbar viel gr\"o\ss{}er ist die Gesamtzahl der Elementarteilchen im Universum, die man auf ungef\"ahr $10^{87}$ sch\"atzt: das ist eine $1$ mit 87 Nullen! Auch diese Zahl k\"onnen wir mit unserer Vorstellungskraft aber immer noch locker \"uberbieten: wir k\"onnen an die Zahl $10^{87}+1$ denken (ein extra Elementarteilchen), oder an $2\times 10^{87}$ (zwei Universen), oder gar an eine 1 mit $10^{87}$ Nullen (so viele Nullen, wie es Elementarteilchen gibt).

Wir k\"onnen uns also beliebig gro\ss{}e Zahlen vorstellen. Anders formuliert gibt es \emph{beliebig viele} nat\"urliche Zahlen 
$$1,\;2,\;3,\;4,\;5,\;6,\;7,\;8,\;9,\;...\;.$$
Man sagt deshalb auch, die nat\"urlichen Zahlen seien \emph{potentiell (also der M\"oglichkeit nach) unendlich}; ihr Vorrat ist unersch\"opflich. 

Tats\"achlich taucht potentielle Unendlichkeit auch in allt\"aglichen Si\-tu\-ationen auf. Wenn man zum Beispiel auf einem Platz steht und sich entscheiden muss, in welche Richtung man weitergeht, dann hat man dabei ganz streng genommen beliebig viele M\"oglichkeiten, n\"amlich so viele wie es Winkel zwischen $0$ und 360 Grad gibt: zum Beispiel $1$ Grad, $\nicefrac{1}{2}$ Grad, $\nicefrac{1}{3}$ Grad, $\nicefrac{1}{4}$ Grad, und beliebig viele Winkel mehr.

\subsection*{Aktuale Unendlichkeit}

Dadurch kommen wir aber noch nicht mit echter Unendlichkeit, der \emph{aktualen Unendlichkeit}, in Ber\"uhrung. Der Vorrat der nat\"urlichen Zahlen ist zwar unbegrenzt, aber ganz konkret ausw\"ahlen, benutzen oder uns vorstellen k\"onnen wir \emph{gleichzeitig} zun\"achst nur \emph{endlich} viele Dinge, seien es auch noch so viele. Oder haben Sie schon einmal unendlich viele Dinge wirklich \emph{auf einmal} gesehen, gesp\"urt, geh\"ort? 

Die gedankliche Zusammenfassung unterschiedlicher Dinge zu einem Gan\-zen nennt man auch \emph{Menge}. Diese (naive) Definition einer Menge geht auf Cantor zur\"uck. Beispielsweise k\"onnen wir die Zahlen $1$, $2$ und $3$ zu einer Menge zusammenfassen, und schreiben daf\"ur $\{1,2,3\}$. Die geschweiften Klammen ``$\{$'' und ``$\}$'' hei\ss{}en auch Mengenklammern, und die Einzelteile, aus der eine Menge besteht, hei\ss{}en \emph{Elemente} der Menge. Damit ist die $2$ ein Element der Menge $\{1,2,3\}$.  

Aktuale Unendlichkeit ist nun etwas unerh\"ort Radikales. Wenn wir die nat\"ur\-lich\-en Zahlen in ihrer Gesamtheit betrachten, das hei\ss{}t in unserer Vorstellung \emph{alle denk\-baren Zahlen auf einmal zusammenfassen}, dann ist die resultierende Menge $\mathbb N$, also 
$$\mathbb N = \{1,\;2,\;3,\;4,\;5,\;6,\;7,\;8,\;9,\;...\}\,,$$
nicht endlich, und in diesem Sinne aktual un-endlich. (Das Symbol ``$\mathbb N$'' ist ein verschn\"orkeltes ``N'' und steht als Abk\"urzung f\"ur ``Nat\"urliche Zahlen''). Eine Menge hei\ss{}t n\"amlich \emph{endlich}, wenn wir alle ihre Elemente z\"ahlen k\"onnen. Bei der Menge $\mathbb N$ k\"onnen wir das ja gerade nicht --  egal wie weit wir z\"ahlen, es bleiben immer noch Elemente \"ubrig.

Ist das zul\"assig? K\"onnen wir das, \emph{d\"urfen} wir einfach unendliche viele Dinge wirklich gleichzeitig betrachten, wenn das keine Entsprechung in der Realit\"at hat? Das ist zun\"achst nicht klar. Wenn \"uberhaupt, dann existiert aktuale Unendlichkeit nur als Konzept in unserer Gedankenwelt.

Die Unterscheidung zwischen potentieller und aktualer Unendlichkeit, welche auf den griechischen Philosophen Aristoteles zur\"uckgeht, ist m\"og\-lich\-er\-weise nicht ganz so leicht zu verstehen. Deshalb nochmal: Es ist ein gewaltiger Unterschied, ob man einerseits eine feste, endliche Menge von Dingen betrachtet und in der Gr\"o\ss{}e dieser endlichen Menge keiner Schranke unterworfen ist; oder ob man andererseits die Endlichkeit unmittelbar aufhebt und in diesem Sinne eine nicht endliche, also unendliche Gesamtheit betrachtet.

\subsection*{Ein kleiner Schritt f\"ur Cantor, ein gro\ss{}er f\"ur die Menschheit}
Bis Cantor hat sich kein Mathematiker oder Philosoph ernsthaft an das Konzept der aktualen Unendlichkeit herangetraut. Einer der bedeutendsten und einfluss\-reich\-sten Mathematiker aller Zeiten, Carl Friedrich Gau\ss{}, schrieb noch 1831:
\begin{quote}
"So protestiere ich zuv\"orderst gegen den Gebrauch einer unendlichen Gr\"o\ss{}e als einer vollendeten [Gr\"o\ss{}e], welches in der Mathematik niemals erlaubt ist."
\end{quote}
Mit einer "vollendeten unendlichen Gr\"o\ss{}e" meint Gau\ss{} gerade die Zusammenfassung unterschiedlicher Dinge zu einem nicht endlichen Ganzen, wie beispielsweise die Menge $\mathbb N$ aller nat\"urlichen Zahlen.

Cantors Arbeiten ab 1870 stie\ss{}en noch bis ins 20. Jahrhundert hinein auf zum Teil heftige Ablehnung unter Mathematikern. Ein bekanntes Zitat dazu, anscheinend ohne eindeutigen Urheber, lautet:
\begin{quote}
"Zuk\"unftige Generationen werden die Mengenlehre als eine Krankheit betrachten, von der man sich erholt hat."
\end{quote}
Heute gilt Cantor als einer der bedeutendsten Mathematiker \"uberhaupt. Von David Hilbert, der die Bedeutung von Cantors Arbeiten fr\"uh erkannt hat, stammt das Zitat (1926):
\begin{quote}
"Aus dem Paradies, das Cantor uns geschaffen, soll uns niemand vertreiben k\"onnen."
\end{quote}

Bis heute gibt es Mathematiker und Philosophen, welche die Idee der aktualen Unendlichkeit ablehnen. Wir m\"ochten uns an dieser Stelle aber nicht weiter mit dieser Diskussion besch\"aftigen, und stellen uns auf den Standpunkt: warum nicht? Schauen wir einfach was passiert.




\subsection*{Ja! -- Aber wie viele?}
Wir sagen also ``Ja!'' zur aktualen Unendlichkeit. In der modernen Mengenlehre wird die Existenz unendlicher Mengen in der Tat als Axiom, das hei\ss{}t als nicht weiter begr\"undete Tatsache angenommen -- in Form des \emph{Unendlichkeits\-axioms}. Darauf kommen wir ganz am Ende nochmal zur\"uck.

Wie k\"onnen wir nun die aktuale Unendlichkeit genauer untersuchen? Wie kann man \"uberhaupt eine \emph{Idee} genauer untersuchen? Da wir hier totales Neuland betreten und uns keine Erfahrung in irgendeiner Weise helfen kann, k\"onnen wir h\"ochstens nach ganz elementaren Eigenschaften fragen. 

Also fragen wir: Besitzt die aktuale Unendlichkeit \emph{\"uberhaupt irgendeine Struktur}? Oder ist Aktualunendlichkeit einfach Aktual\-un\-endlich\-keit, ohne irgendwie geartete Abstufungen, und das ist schon das Ende der ganzen Geschichte? Anders gefragt:
\begin{quote}
\emph{Gibt es verschiedene Arten von Aktualunendlichkeit?}
\end{quote}
Von dieser Frage wollen wir uns im Rest dieses Textes leiten lassen.

\subsection*{Vergleich unendlicher Mengen -- die Grenzen der Intuition}
Eine einfache Beobachtung ist, dass es verschiedene unendliche Mengen gibt. Neben den schon erw\"ahnten nat\"urlichen Zahlen 
$$\mathbb N = \{1,\;2,\;3,\;4,\;5,\;6,\;7,\;8,\;9,\;...\}$$ 
k\"onnen wir zum Beispiel die Menge $\mathbb G$ aller \emph{geraden}, also durch 2 ohne Rest teilbaren Zahlen  betrachten,
$$\mathbb G = \{2,\;4,\;6,\;8,\;10,\;12,\;14,\;16,\;18,\;...\}\,.$$ 
Die Menge $\mathbb G$ erhalten wir, indem wir von $\mathbb N$ jede zweite Zahl weglassen -- nur jede zweite Zahl ist gerade. 

Beide Mengen sind unendlich gro\ss{}, also in einem gewissen Sinne "gleich gro\ss{}". Andererseits hat $\mathbb G$ doch unendlich viele Elemente weniger als $\mathbb N$ und sollte damit  "weniger unendlich" sein als $\mathbb N$. Welche Sichtweise ist nun die richtige? Wie k\"onnen wir \"uberhaupt entscheiden, welche Sichtweise die richtige ist?

Wir wollen dieses Problem  durch ein Beispiel veranschaulichen. Wir stellen uns vor, jemand w\"urde unendlich viel Geld besitzen und m\"usste 50$\%$ Verm\"ogensteuer bezahlen. Das kostet zwar die H\"alfte des Verm\"ogens, aber danach ist ja immer noch unendlich viel Geld \"ubrig. Ist durch die Be\-gleich\-ung der Steuerschuld also \"uberhaupt etwas verloren gegangen? Einerseits unendlich viel, aber andererseits auch wieder nichts. Im Fall eines endlichen Verm\"ogens ist die Sache klar -- wenn man von 100 Euro die H\"alfte abgibt hat man offensichtlich weniger in der Tasche als vorher. 

Anscheinend l\"asst sich diese Argumentation nicht eindeutig auf unend\-liche Mengen \"ubertragen, unsere Intuition verl\"asst uns an dieser Stelle. Wir k\"onnen ja gerade nicht die Elemente unendlicher Mengen z\"ahlen und die Ergebnisse vergleichen!


\subsection*{Gleichm\"achtigkeit}
Wie k\"onnen wir also die "Gr\"o\ss{}e" unendlicher Mengen miteinander ver\-gleichen? Cantors Idee zur L\"osung dieses Problems ist so einfach wie genial. Wir m\"ochten ja gar nicht wissen, "wie viele" Elemente eine unendliche Menge besitzt, denn diese Information gibt es schlicht und einfach nicht. Wir wollen unendliche Mengen, zum Beispiel $\mathbb N$ und $\mathbb G$, ja nur \emph{vergleichen}. 

Cantors Idee dazu basiert auf der Tatsache, dass es endlichen Mengen au\ss{}er dem Z\"ahlen-und-Ergebnisse-Vergleichen noch eine andere Methode zum Gr\"o\ss{}envergleich gibt, die nicht nach der Anzahl der Elemente fragt: \emph{Paarbildung ohne Rest}.

Stellen wir uns vor, wir m\"ochten wissen, ob beim Wiener Opernball w\"ahrend einer Pause mehr Frauen oder mehr M\"anner auf der Tanzfl\"ache sind. Man kennt die Bilder -- schwarze Fracks und wei\ss{}e Ballkleider wuseln v\"ollig un\"ubersichtlich durcheinander. Nun k\"onnten wir erst die Frauen durchz\"ahlen, dann die M\"anner, und am Ende die Ergebnisse vergleichen. Das liefert uns aber viel mehr als wir eigentlich wissen wollten, n\"amlich sogar die genaue Anzahl der Frauen und M\"anner. 

Eine zweite, einfachere M\"oglichkeit ist, \"uber die Lautsprecher durchzu\-sagen, jede Frau solle sich zu einem Mann gesellen. Sind danach noch Frauen \"ubrig, dann wissen wir sicher, dass mehr Frauen auf der Tanzfl\"ache sind. Sind noch M\"anner \"ubrig, dann sind es mehr M\"anner. Bleibt niemand alleine stehen, d.h. es gelang eine Paarbildung ohne Rest, dann sind gleich viele Frauen und M\"anner auf der Tanzfl\"ache.

Schon Kinder wenden das Prinzip der Paarbildung ohne Rest zum Gr\"o\ss{}en\-vergleich erfolgreich an. Eine Schachtel Pralinen werden sie niemals so aufteilen, dass sie die Pralinen z\"ahlen und das Ergebnis halbieren, sondern eher nach dem Prinzip: eins f\"ur dich, eins f\"ur mich, eins f\"ur dich, und so weiter.

Diese M\"oglichkeit zum Gr\"o\ss{}envergleich bei endlichen Mengen erweiterte Cantor nun konsequent durch folgende \emph{Definition} auf beliebige, auch unendliche Mengen. Um aber nicht mit anderen Bedeutungen von "Gr\"o\ss{}e" im Zusammenhang mit Mengen in Konflikt zu geraten, m\"ochten wir wie Cantor lieber den Begriff "M\"achtigkeit" benutzen. 

\begin{quote}
Zwei beliebige Mengen $A$ und $B$ hei\ss{}en \emph{gleichm\"achtig}, wenn man Paare aus ihren Elementen bilden kann, so dass kein Element aus $A$ und $B$ \"ubrig bleibt. Symbolisch schreibt man daf\"ur $|A| = |B|$.
\end{quote}
In dieser Definition sind $A$ und $B$ symbolische Platzhalter f\"ur beliebige Mengen unserer Wahl. Zum Beispiel k\"onnten wir $A = \mathbb N$ und $B=\mathbb G$ setzen um die M\"achtigkeit der nat\"urlichen und der geraden Zahlen zu vergleichen.

Wir haben oben gesehen, dass bei endlichen Mengen $A$ und $B$ Cantors Definition von Gleichm\"achtigkeit zu dem zur\"uckf\"uhrt, was man \"ublichweise unter zwei gleichgro\ss{}en Mengen versteht. Da sie ferner auf beliebige Mengen anwendbar ist, \emph{glauben} wir, dass diese Definition sinnvoll ist. 

Jedoch ist an dieser Stelle noch gar nicht klar, ob die Definition auch "gut" ist, das hei\ss{}t, ob der Begriff der Gleichm\"achtigkeit uns wirklich etwas Neues \"uber die aktuale Unendlichkeit lehrt. Es k\"onnte ja sein, dass sich herausstellt, dass alle unendlichen Mengen gleichm\"achtig sind. Dann w\"urde dieser Begriff keine Strukturen innerhalb der Aktualunendlichkeit beschreiben. Vorab darf aber verraten werden, dass Cantors De\-finition uns tats\"achlich neue Einsichten \"uber die Unendlichkeit liefert.

Bevor wir uns aber im n\"achsten Abschnitt an unendliche Mengen wagen, wollen wir die Paarbildung ohne Rest zun\"achst anhand endlicher Mengen veranschaulichen. Die Mengen $A = \{1,2,3\}$ und $B= \{a,b,c\}$ haben jeweils drei Elemente, sind also gleichgro\ss{}. Eine Paarbildung ohne Rest sieht zum Beispiel so aus:
$$1\leftrightarrow a,\qquad 2\leftrightarrow b, \qquad 3\leftrightarrow c.$$
Genauso ist auch
$$1\leftrightarrow b,\qquad 2\leftrightarrow c, \qquad 3\leftrightarrow a,$$
eine Paarbildung ohne Rest. Dies zeigt an einem ganz einfachen Beispiel, wie man die Gleichm\"achtigkeit zweier Mengen durch Paarbildung ohne Rest nachweist.

\subsection*{Gr\"o\ss{}ere und kleinere M\"achtigkeiten}

Da wir nun festgelegt haben, wann wir zwei Mengen gleichm\"achtig nennen wollen, k\"onnen wir auf Basis derselben Idee definieren, wann eine beliebige Menge "m\"achtiger" (also gr\"o\ss{}er) als eine andere sein soll:

\begin{quote}
Die Menge $B$ hei\ss{}t \emph{m\"achtiger} als die Menge $A$, wenn $A$ gleich\-m\"achtig zu einer Teilmenge von $B$ ist, aber $A$ und $B$ nicht gleich\-m\"achtig sind. Symbolisch schreibt man daf\"ur $|A| < |B|$.
\end{quote}

Halten wir kurz inne und schauen uns die Definition genauer an. Dazu, dass $B$ m\"achtiger als $A$ ist, fordern wir zwei Dinge. Zun\"achst sollen $A$ und $B$ \"uberhaupt einmal nicht gleichm\"achtig sein, es soll also keine Paarbildung ohne Rest geben. Klar, dass dies erf\"ullt sein sollte, sonst k\"onnen wir ja nur schlecht $B$ m\"achtiger als $A$ nennen. Nun stellt sich die Frage, ob $A$ oder $B$ die "gr\"o\ss{}ere" Menge ist. Daf\"ur ist der andere Punkt da: Wenn wir von $B$ etwas Passendes wegnehmen, also nur eine bestimmte Teilmenge von $B$ betrachten, dann sollen wir eine zu $A$ gleichm\"achtige Menge erhalten. In diesem Fall k\"onnen wir $B$ mit Fug und Recht als m\"achtiger, also "gr\"o\ss{}er", als $A$ bezeichnen.

Veranschaulichen wir dies anhand der Mengen $A = \{1,2,3\}$ und $B = \{a, b, c, d\}$. Sofort sehen wir, dass $B$ im \"ublichen Sinne gr\"o\ss{}er ist als $A$, da es ein Element mehr als $A$ enth\"alt. Aber darum geht es uns bei diesem einfachen Beispiel ja wieder  nicht. Wir wollen uns nur auf Cantors Definitionen verlassen, denn f\"ur die F\"alle, die uns wirklich interessieren, also unendliche Mengen, haben wir nur diese zur Verf\"ugung!

Wir wollen im Cantorschen Sinne zeigen, dass $B$ m\"achtiger als $A$ ist. Wie wir oben gesehen haben, ist $A$ gleichm\"achtig zur Teilmenge $\{a,b,c\}$ von $B$. Das w\"are also erledigt. Wir m\"ussen nun noch nachweisen, dass $A$ und $B$ nicht gleichm\"achtig sind, also \emph{keine} Paarbildung ohne Rest m\"oglich ist. In diesem einfachen Beispiel ist das klar, oder?  Nun ja, erstens m\"ochten wir uns als Mathematiker alles ganz genau \"uberlegen und  bei allem 100$\%$ sicher sein (Stichworte Transparenz und Rigorosit\"at). Zweitens wird bei unendlichen Mengen \"uberhaupt nichts mehr klar sein, so dass wir gut daran tun, jetzt schon damit anfangen zu \"uberlegen, wie wir eine Paarbildung ohne Rest auf jeden Fall aussschlie\ss{}en k\"onnen.

Eine todsichere Methode ist es, uns die M\"uhe zu machen, alle m\"oglichen Paarbildungen durchzugehen (insgesamt $24 = 4\times 3\times 2$ St\"uck) und zu \"uberpr\"ufen, ob immer ein Rest bleibt. Ein Beispiel ist
$$1\leftrightarrow c, \qquad 2\leftrightarrow d,\qquad 3\leftrightarrow a, \qquad ?\leftrightarrow b.$$
Also bleibt das Element $b$ von $B$ ohne Partner aus $A$. Nun mag man sich davon \"uberzeugen, dass es bei jeder der 23 anderen Paarbildungen genauso l\"auft. Damit folgern wir auch aus Cantors Definition, dass $B$ m\"achtiger ist als $A$.

Das Beispiel zeigt, dass der Aufwand Cantors Definition von "m\"achtiger" nachzupr\"ufen, sehr gro\ss{} sein kann, was bei einfachen Beispielen zurecht unverh\"altnism\"a\ss{}ig erscheinen mag. Im Fall von endlichen Mengen kann man sich \"uberlegen, dass eine Paarbildung wie die andere ist -- entweder alle Paarbildungen sind ohne Rest, oder es bleibt bei jeder Paarbildung mindestens ein Element \"ubrig, und zwar immer von der selben Menge. Das k\"urzt das ganze etwas ab.

Um endliche Mengen geht es uns aber gar nicht. Viel wichtiger ist, dass Cantors Definition im Fall unendlicher Meng\-en \emph{die einzige M\"oglichkeit} zu sein scheint, \"uberhaupt einen sinnvollen Gr\"o\ss{}\-en\-vergleich durchzuf\"uhren.

Wir sollten uns auch nochmal klar machen, dass wir an dieser Stelle nur das \emph{Instrument} besitzen, Gr\"o\ss{}en\-ab\-stufungen innerhalb der Aktualunendlich\-keit und damit unterschiedliche Arten von Unendlichkeit zu finden. Noch ist unklar, ob es \"uberhaupt unendliche Mengen mit verschiedenen M\"achtigkeiten gibt, oder ob einfach alle unendlichen Mengen gleichm\"achtig sind.

\subsection*{Die Gleichm\"achtigkeit von $\mathbb N$ und $\mathbb G$}
Wir wollen nun im Folgenden mit Hilfe des Instruments, das uns Cantor durch seine Definitionen gegeben hat, die M\"achtigkeit unendlicher Mengen miteinander vergleichen.

Wir beginnen mit $\mathbb N$ und $\mathbb G$. Tats\"achlich finden wir eine Paarbildung ohne Rest, basierend auf folgender Idee: die 2 ist die \emph{erste} gerade Zahl, die 4 ist die \emph{zweite} gerade Zahl, die 6 die \emph{dritte}, und so weiter. In einem Schaubild sieht unsere Paarbildung  so aus:
$$\begin{array}{rcccccccccc}
\mathbb N :\;\;& 1 & 2 & 3 & 4 & 5 & 6 & 7 & 8 &9  &\cdots\\
& \updownarrow & \updownarrow &\updownarrow &\updownarrow &\updownarrow &\updownarrow &\updownarrow &\updownarrow &\updownarrow &\\
\mathbb G :\;\;& 2 & 4 & 6 & 8 & 10 & 12 & 14 & 16 & 18  &\cdots
\end{array}
$$
Kritische Leser werden bem\"angeln, dass die P\"unktchen "$\cdots$" etwas ungenau sind. In der Tat. Anders als bei endlichen Mengen k\"onnen wir nicht alle Paare aus $\mathbb N$ und $\mathbb G$ auf einmal hinschreiben, es sind ja unendlich viele. Die \emph{Vorschrift} zur Paarbildung ist aber endlich: jede nat\"urliche Zahl wird mit ihrem Doppelten gepaart. Damit hat jede Zahl ihren eindeutigen ge\-raden Partner, und wir erwischen alle geraden Zahlen (hier finden wir einen zugeh\"origen Partner aus $\mathbb N$ jeweils durch Halbierung). Das reicht uns als ordentliche Paarbildung wie in Cantors Definition verlangt.

Es folgt, dass $\mathbb N$ und $\mathbb G$ gleichm\"achtig  sind,
$$|\mathbb N| = |\mathbb G|.$$
Wenn wir Cantors Definition f\"ur sinnvoll und gut halten (siehe oben), m\"ussen wir einsehen,  dass es "genauso viele" nat\"urliche wie gerade Zahlen gibt! Und das, obwohl $\mathbb G$ eine echte Teilmenge von $\mathbb N$ ist! Wir halten fest:

\begin{quote}
\emph{Durch Wegnehmen von unendlich vielen Elementen wird eine unendliche Menge nicht unbedingt weniger m\"achtig.}
\end{quote}

Dies ist bei endlichen Mengen nicht m\"oglich. Durch Zulassen unendlicher Mengen sto\ss{}en wir also sofort auf ganz neue, unbekannte Ph\"ano\-mene. Wir nehmen jedes zweite Element von $\mathbb N$ weg, und trotzdem \"andern wir nichts an der M\"achtigkeit bzw. "Gr\"o\ss{}e"! Nat\"urlich hei\ss{}t das nicht, dass \emph{ganz egal} wie viele Elemente wir von $\mathbb N$ wegnehmen sich die M\"achtigkeit nicht \"andert. Nehmen wir von $\mathbb N$ alles weg au\ss{}er der $1$, dann haben wir die M\"achtigkeit  wirklich verkleinert. Aber wie ist es, wenn wir noch unendlich viele Elemente \"ubrig lassen? Dazu gleich mehr.

Zun\"achst noch eine Bemerkung. Ein weiteres neues Ph\"anomen bei unendlichen Mengen ist, dass es noch andere Paarbildungen gibt, bei denen sogar unendlich viele Elemente zum Beispiel von $\mathbb N$ \"ubrig bleiben, obwohl wir schon eine erfolgreiche Paarbildung ohne Rest gefunden haben. Beispiels\-weise lassen wir durch die Paarbildung
$$\begin{array}{rcccccccccc}
\mathbb N :\;\;& 1 & 2 & 3 & 4 & 5 & 6 & 7 & 8 &9  &\cdots\\
&  & \updownarrow & &\updownarrow & &\updownarrow & &\updownarrow & &\\
\mathbb G :\;\;&  & 2 &  & 4 &  & 6 & & 8 &   &\cdots
\end{array}
$$
genau die ungeraden Zahlen  aus, obwohl wir alle Elemente von $\mathbb G$ ber\"uck\-sichtigt haben. Dies widerspricht aber nicht der Definition von Gleich\-m\"achtig\-keit, denn bei dieser wird ja nur nach \emph{mindestens einer} erfolgreichen Paarbildung ohne Rest gefragt. Ob es noch andere Paarbildungen \emph{mit Rest} gibt, spielt keine Rolle.

\subsection*{Die M\"achtigkeit von $\mathbb N$ als kleinste unendliche M\"achtigkeit}

Erstaunlicherweise gibt es also gleich viele nat\"urliche wie gerade Zahlen, obwohl nur jede zweite nat\"urliche Zahl gerade ist. Wir wollen dies konsequent weiterdenken und fragen, was passiert, wenn wir einfach noch mehr von $\mathbb N$ wegnehmen. Erhalten wir dann irgendwann eine unendliche Teilmenge von $\mathbb N$, die weniger m\"achtig ist als $\mathbb N$? Noch allgemeiner gefragt: Gibt es \"uberhaupt eine unendliche Menge, die weniger m\"achtig ist als die nat\"urlichen Zahlen?

Dem ist in der Tat nicht so. Um dies einzu\-sehen, nehmen wir uns irgendeine unendliche Menge $A$ her, die nicht gleichm\"achtig zu $\mathbb N$ ist. Nun betrachten wir folgende Paarbildung: als Partner der $1$ w\"ahlen wir irgend\-ein Element von $A$, als Partner der $2$ irgendein weiteres Element von $A$, als Partner der $3$ noch irgendein weiteres Element von $A$, und so weiter. Dadurch bekommt jedes Element von $\mathbb N$ einen Partner ab. Warum? Wenn es z.B. ab der Zahl 1000 keine weiteren Elemente in $A$ g\"abe, dann h\"atte $A$ ja nur 1000 Elemente gehabt, w\"are also nicht unendlich. 

Zwisch\-en $\mathbb N$ und den Elementen aus $A$, die wir in diesem Prozess ausgew\"ahlt haben, erhalten wir also eine Paarbildung ohne Rest. Damit sind $\mathbb N$ und diese Teilmenge von $A$ gleichm\"achtig. Insgesamt sind $\mathbb N$ und $A$ also nicht gleichm\"achtig und $\mathbb N$ ist gleichm\"achtig zu einer Teilmenge von $A$. Damit ist $A$ laut Cantors Definition m\"achtiger als $\mathbb N$. Wir folgern: 
\begin{quote}
\emph{Jede unendliche Menge ist mindestens so m\"achtig wie $\mathbb N$.}
\end{quote}
\emph{Wenn} es also \"uberhaupt unterschiedliche unendliche M\"achtigkeiten gibt, dann gibt es auch eine kleinste, n\"amlich die M\"achtigkeit der nat\"urlichen Zahlen. Diese Erkenntnis ist fundamental. Die M\"achtigkeit der nat\"urlichen Zahlen ist also etwas ganz Besonderes, weshalb wir einen eigenen Begriff daf\"ur haben m\"ochten.

\subsection*{Abz\"ahlbarkeit}

Wie schon erw\"ahnt kann man die oben gezeigte Paarbildungsvorschrift ohne Rest zwischen $\mathbb G$ und $\mathbb N$ auch als ein "unendliches Abz\"ahlen" der Elemente von $\mathbb G$ interpretieren. Den Partner der  $1$ fassen wir als die \emph{erste} gerade Zahl auf, den Partner der $2$ als die \emph{zweite} gerade Zahl, den Partner der $3$ als die \emph{dritte}, und so weiter. Anders gesagt kann man die Elemente einer Menge, die gleich\-m\"achtig zu $\mathbb N$ ist, in eine unendliche Liste schreiben: den Partner der $1$ als erstes Listenelement, den Partner der $2$ als zweites, und so weiter. 

Weil die M\"achtigkeit der nat\"urlichen Zahlen so fundamental ist, m\"ochten wir deshalb im Folgenden eine unendliche Menge, die gleichm\"achtig zu $\mathbb N$ ist, auch \emph{abz\"ahlbar unendlich} oder auch \emph{auflistbar unendlich} nennen. Meist sagen wir auch einfach nur \emph{abz\"ahlbar}. Zum Beispiel ist $\mathbb G$ ab\-z\"ahl\-bar, und gleich werden wir noch mehr abz\"ahl\-bare Mengen kennenlernen.

Der Begriff der Abz\"ahlbarkeit sollte nicht mit ``Z\"ahlbarkeit'' verwechselt werden. Letzteres bedeutet n\"amlich, dass man die Elemente einer Menge z\"ahlen kann, dass sie also nur endlich viele Elemente hat. Im Gegensatz dazu entspricht Abz\"ahlbarkeit wie schon beschrieben einem ``unendlichen Abz\"ahlen''.

Eine Menge, die m\"achtiger ist als $\mathbb N$, k\"onnen wir nicht abz\"ahlen bzw. auflisten. Deshalb nennen wir solche Menge  \emph{\"uberabz\"ahlbar unendlich} oder einfach \emph{\"uberabz\"ahlbar}.

\subsection*{Weitere abz\"ahlbare Mengen}

Nachdem wir gesehen haben, dass die M\"achtigkeit von $\mathbb N$ die kleinste unendliche M\"achtigkeit ist, wollen wir nun Elemente zu $\mathbb N$ hinzuf\"ugen und \"uberpr\"ufen, ob sich dadurch die M\"achtigkeit irgendwann ver\"andert, ob wir also irgendwann eine \"uberabz\"ahlbar unendliche Menge erhalten.

Beginnen wir ganz bescheiden mit einem zus\"atzlichen Element und betrachten
$$\mathbb N_0 = \{0,\;1,\;2,\;3,\;4,\;5,\;6,\;7,\;8,\;...\},$$
also die nat\"urlichen Zahlen $\mathbb N$ zusammen mit der Null.  Eine Paarbildung ohne Rest sieht so aus:
$$\begin{array}{rcccccccccc}
\mathbb N :\;\;& 1 & 2 & 3 & 4 & 5 & 6 & 7 & 8 &9  &\cdots\\
& \updownarrow & \updownarrow &\updownarrow &\updownarrow &\updownarrow &\updownarrow &\updownarrow &\updownarrow &\updownarrow &\\
\mathbb N_0 :\;\;& 0 & 1 & 2 & 3 & 4 & 5 & 6 & 7 & 8  &\cdots
\end{array}
$$
Die Paarbildungvorschrift lautet also: Ziehe immer 1 ab. Damit sind $\mathbb N$ und $\mathbb N_0$ gleichm\"achtig, $|\mathbb N| = |\mathbb N_0|$, das hei\ss{}t auch $\mathbb N_0$ ist abz\"ahlbar unendlich. 

Genauso \"uberlegt man sich, dass das Hinzuf\"ugen von 19 oder einer Million neuer Elemente die M\"achtigkeit von $\mathbb N$ nicht \"andert. Die Paarbildungsvorschrift lautet hier entsprechend: ziehe immer 19 oder eine Million ab. Allgemeiner \"andert sich die M\"achtigkeit durch Hinzuf\"ugen von endlich vielen Elementen nicht. 

Im n\"achsten Schritt nehmen wir unendlich viele Elemente zu $\mathbb N$ hinzu und betrachten die \emph{ganzen Zahlen} $\mathbb Z$,
$$\mathbb Z = \{..., \;-4,\;-3,\;-2,\;-1,\;0,\;1,\;2,\;3,\;4,\;...\},$$
die neben den nat\"urlichen Zahlen auch alle negativen Zahlen umfassen. Nun k\"onnten wir vielleicht eine \"uberabz\"ahlbare Menge gefunden haben, denn immerhin ist $\mathbb Z$ "nach links und rechts schrankenlos", w\"ahrend $\mathbb G$ und $\mathbb N_0$ nur in einer Richtung schrankenlos sind. Darauf kommt es aber nicht an. Indem wir $\mathbb Z$ einfach anders aufschreiben, n\"amlich als
$$\mathbb Z = \{0,\;1,\;-1,\;2,\;-2,\;3,\;-3,\;4,\;-4,\;...\},$$
wird wieder klar, wie eine Paarbildung mit $\mathbb N$ ohne Rest geht:
$$\begin{array}{rcccccccccc}
\mathbb N :\;\;& 1 & 2 & 3 & 4 & 5 & 6 & 7 & 8 &9  &\cdots\\
& \updownarrow & \updownarrow &\updownarrow &\updownarrow &\updownarrow &\updownarrow &\updownarrow &\updownarrow &\updownarrow & \\
\mathbb Z :\;\;& 0 & 1 & -1 & 2 & -2 & 3 & -3 & 4 & -4  &\cdots
\end{array}
$$
Die dahintersteckende endliche Paarbildungsvorschrift unterscheidet also zwi\-schen geraden und ungeraden Elementen von $\mathbb N$. Eine ge\-rade Zahl wird halbiert, so kommen wir zum Beispiel von der  $4$ zur $2$. Von einer ungeraden Zahl zieht man 1 ab, halbiert und setzt ein Minus davor. Dies ergibt zum Beispiel die $-4$ als Partner der $9$, denn $-4 = -(9-1)/2$.  Um sicher zu sein, dass wir durch diese Art der Paarbildung wirklich jedes Element von $\mathbb Z$ erwischen, geben wir einfach den Partner aus $\mathbb N$ einer ganzen Zahl konkret an: eine positive Zahl wird verdoppelt; bei einer negativen Zahl wird das Minus weggelassen, dann verdoppelt und schlie\ss{}lich eins dazuaddiert.

Also ist auch $\mathbb Z$ abz\"ahlbar, es gibt genauso viele ganze wie nat\"urliche Zahlen! Nachdem wir $\mathbb Z$ anderes aufgeschrieben haben, war das auch schon so zu erwarten. Der Schritt von $\mathbb Z$ zu $\mathbb N$ ist analog zum Schritt von $\mathbb N$ nach $\mathbb G$, n\"amlich Weglassen jedes zweiten Elements, und die Gleichm\"achtigkeit von $\mathbb N$ und $\mathbb G$ hatten wir ja schon nachgewiesen. Daraus folgt:

\begin{quote}
\emph{Durch Hinzuf\"ugen von unendlich vielen Elementen wird eine unendliche Menge nicht unbedingt m\"achtiger.}
\end{quote}

Schlie\ss{}lich wollen wir an dieser Stelle festhalten, dass alle unendlichen Mengen, die wir bis hierhin untersucht haben, abz\"ahlbar sind:
$$|\mathbb N| = |\mathbb G| = |\mathbb N_0| = |\mathbb Z|.$$

\subsection*{Ein Ausflug in Hilberts Hotel}
Wir m\"ochten nun in einem kurzen Ausflug die Eigenschaften unendlicher Mengen, die wir bis jetzt kennengelernt haben, etwas mehr veranschaulichen. Von David Hilbert stammt dazu ein sehr sch\"ones Gedankenexperiment, welches heute als \emph{Hilberts Hotel} bekannt ist.

Stellen wir uns vor, es g\"abe ein Hotel mit abz\"ahlbar unendlich vielen Zimmern, also so viele Zimmer, wie es nat\"urliche Zahlen gibt. Stellen wir uns weiterhin vor, alle Zimmer des Hotels w\"aren belegt. Was, wenn nun noch ein Gast auftaucht und nach einem freien Zimmer fragt? M\"ussen wir ihm dann sagen, es w\"are leider keines mehr frei und ihn an ein anderes Hotel verweisen? Auf den ersten Blick ja, denn es sind ja alle Zimmer belegt. 

Andererseits haben wir es ja nicht mit einem normalen, endlichen Hotel zu tun. Im letzten Abschnitt haben wir gesehen, dass wir durch Hinzuf\"ugen eines Elements zu $\mathbb N$ die M\"achtigkeit nicht \"andern. Bei der zugeh\"origen Paarbildungvorschrift sind alle nat\"urlichen Zahlen einfach eins r\"uberger\"uckt, und schon war Platz f\"ur ein weiteres Element, die Null. 

Genauso k\"onnen wir auch in unserem unendlichen Hotel noch Platz f\"ur einen neuen Gast schaffen: wir bitten den Gast aus Zimmer 1 in Zimmer 2 umzuziehen, den Gast aus Zimmer 2 in Zimmer 3, den Gast aus Zimmer 3 in Zimmer 4, und so weiter. Dann sind alle fr\"uheren G\"aste wieder in einem Zimmer -- und Zimmer 1 ist frei f\"ur den neuen Gast!

\"Ahnlich k\"onnen wir verfahren, wenn 346 oder eine Milliarde neue G\"aste in Hilberts Hotel neu unterzubringen sind.

Und wenn gleich\-zeitig \emph{abz\"ahlbar unendlich viele neue G\"aste} nach einem Zimmer fragen? Auch dann m\"ussen wir niemanden abweisen. Inspi\-riert von unserer Paarbildungvorschrift ohne Rest zwischen $\mathbb N$ und $\mathbb G$ bitten wir jeden der aktuellen G\"aste in das Zimmer mit der \emph{doppelten Zimmernummer} umzuziehen. Nach Ende dieser Umziehaktion sind nur noch die Zimmer mit gerader Zimmernummer besetzt. Also sind alle Zimmer mit ungerader Zimmernummer freigeworden und wir haben Platz f\"ur so viele G\"aste, wie es ungerade Zahlen gibt. Da es genauso viele ungerade wie nat\"urliche Zahlen gibt (wie lautet hier eine Paarbildungsvorschrift ohne Rest?) ist nun Platz f\"ur abz\"ahlbar unendlich viele G\"aste!

Kurz zusammengefasst: Ein Hotel mit (abz\"ahlbar) unendlich vielen Zimmern kann niemals ausgebucht sein, selbst bei (abz\"ahlbar) unendlich vielen neuen G\"asten.

Nat\"urlich l\"asst sich das nicht mit unserer ``end\-lichen" Anschauung ver\-einbaren. Ist aber kein Problem: es gibt in der Realit\"at einfach keine Hotels mit unendlich vielen Zimmern. Denn wir d\"urfen nicht vergessen, von wo wir gestartet sind, n\"amlich von der Grund\-annahme, dem Axiom, dass es aktualunendliche Mengen wie zum Beispiel die Gesamtheit der nat\"urlichen Zahlen gibt. Schon diese Annahme existiert nur als Konzept. Des\-halb ist es nicht verwunderlich, dass diese Annahme Konsequenzen nach sich zieht, die auch nur in der Gedankenwelt existieren.

Um die Geschichte noch auf die Spitze zu treiben: Einige G\"aste werden wegen der vielen Umzieherei zwar meckern und im Internet schlechte Bewertungen abgeben. Das Preis-Leistungs-Verh\"altnis von Hilberts Hotel ist aber unschlagbar, denn wir k\"onnen die Zimmer sogar kostenlos anbieten. Wie geht das? Wenn wir Geld ben\"otigen, bitten wir einfach den Gast in Zimmer 1 darum. Dieser fragt den Gast aus Zimmer 2 nach demselben Betrag, dieser wiederum den Gast aus Zimmer 3, und so weiter. So k\"onnen wir gut vom Hotel leben und alle unsere Rechnungen bezahlen, und kein Gast ist auch nur um einen Kreuzer \"armer geworden...

\subsection*{Die Abz\"ahlbarkeit der rationalen Zahlen -- \\Cantors erstes Diagonalargument}
Wir setzen unsere Suche nach einer \"uberabz\"ahlbaren Menge fort. Neben den negativen Zahlen nehmen wir jetzt noch alle Br\"uche zu $\mathbb N$ hinzu und erhalten dadurch die \emph{rationalen Zahlen} $\mathbb Q$,
$$\mathbb Q = \left \{\nicefrac{p}{q}\;:\; \text{$p$ und $q$ aus $\mathbb Z$ mit $q$ ungleich $0$}\right\}.$$
Beispielsweise sind $\nicefrac{1}{3}$, $-\nicefrac{7}{28}$ und $\nicefrac{1000}{9}$ rationale Zahlen. Ferner sind auch alle ganzen Zahlen rational, da z.B. $5 = \nicefrac{5}{1}$ und $-2 = \nicefrac{-2}{1}$ gilt.

\"Ubrigens benutzt man hier ein verschn\"orkeltes "Q" als Abk\"urzung f\"ur die rationalen Zahlen, weil das verschn\"orkelte "R", also $\mathbb R$, schon f\"ur die noch wichtigeren reellen Zahlen reserviert ist, auf die wir bald zur\"uckkommen werden.

Nun scheinen wir wirklich eine neue Art von Unendlichkeit gefunden zu haben. Zwischen den Elementen der Mengen, die wir bis jetzt betrachtet haben, klafften jeweils L\"ucken: so ist zwischen $4$ und $5$ keine weitere nat\"urliche Zahl. Bei den rationalen Zahlen ist das ganz anders: zwischen je zwei rationalen Zahlen liegt noch eine weitere rationale Zahl, es gibt keine Nachfolger! Zum Beispiel liegt $\nicefrac{1}{2}$ zwischen $0$ und $1$, es liegt $\nicefrac{1}{3}$ zwischen $0$ und $\nicefrac{1}{2}$, etc. Konsequenterweise liegen zwischen je zwei rationalen Zahlen sogar nochmal unendlich viele weitere rationale Zahlen! $\mathbb Q$ ist also ziemlich dicht gepackt. Es sieht so aus, als g\"abe es sehr viel ``mehr'' Br\"uche als nur nat\"urliche Zahlen.

Auch hier t\"auscht uns die Anschauung. Cantor hat gezeigt, dass die rationalen Zahlen abz\"ahlbar sind. 

Welche Paarbildungsvorschrift ohne Rest zwischen $\mathbb N$ und $\mathbb Q$ hat Cantor  gefunden? Der Einfachheit halber wollen wir nur eine solche zwischen $\mathbb N$ und den positiven rationalen Zahlen angeben und \"uberlassen es dem interessierten Leser, sich zu \"uberlegen, wie daraus schon $|\mathbb Q| = |\mathbb N|$ folgt.

Cantors geniale Idee, wie man die (positiven) rationalen Zahlen abz\"ahlen kann, basiert  darauf, das folgende unendliche rechteckige Schema von Br\"uchen zu betrachten, 
$$
\begin{array}{cccccc}
\nicefrac11  &   \nicefrac12  & \nicefrac13  &\nicefrac14      & \cdots \\\\
\nicefrac21  & \nicefrac22  & \nicefrac23  &\nicefrac24      & \cdots \\\\
\nicefrac31  & \nicefrac32  & \nicefrac33  &\nicefrac34     & \cdots \\\\
\nicefrac41  & \nicefrac42  & \nicefrac43  &\nicefrac44 & \cdots \\
\vdots  & \vdots  & \vdots   & \vdots &  &
\end{array}
$$
In diesem Schema taucht jede positive rationale Zahl mindestens einmal auf. Nun erhalten wir nat\"urlich keine Paarbildung ohne Rest, wenn wir links oben bei $\nicefrac11$ anfangen und in die vertikale oder horizontale Richtung durchz\"ahlen: damit kommen wir nicht mal in die N\"ahe der meisten Br\"uche. Viel besser ist es, \emph{diagonal} abzuz\"ahlen -- das ist gerade \emph{Cantors erstes Diagonalargument}. Und das geht so:
$$
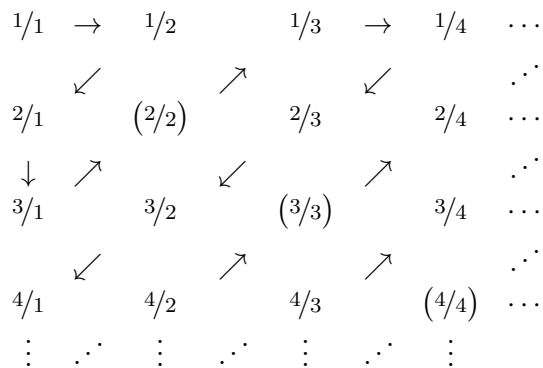
\begin{array}{ccccccccc}
\nicefrac11  &\rightarrow &   \nicefrac12  & &\nicefrac13  & \rightarrow & \nicefrac14      & \cdots \\
& \swarrow & & \nearrow & & \swarrow & & \iddots\\
\nicefrac21  & &\big (\nicefrac22\big)  & &\nicefrac23  & &\nicefrac24      & \cdots \\
\downarrow & \nearrow & & \swarrow & & \nearrow & & \iddots \\
\nicefrac31  & &\nicefrac32  & &\big(\nicefrac33\big)  & &\nicefrac34     & \cdots \\
& \swarrow & & \nearrow  &&\nearrow && \iddots\\
\nicefrac41  & &\nicefrac42  &  &\nicefrac43  && \big(\nicefrac44\big) & \cdots \\
\vdots  &\iddots & \vdots  & \iddots&\vdots   & \iddots & \vdots &  &
\end{array}
$$
Wir beginnen also ganz links oben, der Partner der $1$ ist $\nicefrac11 = 1$. Die Partner von $2$ und $3$ sind die beiden Eintr\"age auf der zweiten Diagonalen, n\"amlich $\nicefrac12$ und $\nicefrac21$. Die Partner von $4$ und $5$ sind die beiden Eintr\"age $\nicefrac31$ und $\nicefrac13$ auf der dritten Diagonalen (hier m\"ussen wir $\nicefrac22 = 1$ nicht beachten, da wir ja schon einen Partner f\"ur die $1$ haben). Die Partner von $6$, $7$, $8$ und $9$ sind die Eintr\"age $\nicefrac14$, $\nicefrac23$, $\nicefrac32$ und $\nicefrac41$ der vierten Diagonalen. Und so weiter. Auf diese Art und Weise schl\"angeln wir uns auf den Diagonalen entlang durch das gesamte Zahlenschema und finden f\"ur jede nat\"urliche Zahl einen Partner aus $\mathbb Q$ (und umgekehrt).

Damit gibt es genauso viele nat\"urliche wie rationale Zahlen, $$|\mathbb N| = |\mathbb Q|.$$


\subsection*{Das Kontinuum -- die Menge der reellen Zahlen}
Nun geben wir endg\"ultig jegliche Zur\"uckhaltung auf und betrachten \emph{alle Dezimalzahlen auf einmal}. Eine Dezimalzahl ist anschaulich eine "Kommazahl mit beliebiger Zahlenfolge hinterm Komma". Beispiele daf\"ur sind alle rationalen Zahlen, zum Beispiel
$$\nicefrac13 = 0,33333...\;, \qquad 5 = 5,00000...\;, \qquad -\nicefrac65 = -1,20000...\;;$$
ferner alle Wurzeln, wie
$$\sqrt{2} = 1,41421...\;;$$
und schlie\ss{}lich alle \"ubrigen Dezimalzahlen, die  sogenannten "transzendenten Zahlen", wie die Kreiszahl $\pi$ und die Eulersche Zahl $e$,
$$\pi = 3,14159...\;,\qquad e = 2,71828...\;.$$
Die resultierende Menge aller Dezimalzahlen nennt man \emph{reelle Zahlen} und k\"urzt sie mit $\mathbb R$ ab. Wie bei $\mathbb N$, $\mathbb Z$ und $\mathbb Q$ kann man nun $\mathbb R$ mit Mengenklammern schreiben als $\mathbb R = \{ ... \}$. Dies sieht aber recht kompliziert aus und ist im Folgenden nicht so wichtig. Merken sollten wir uns einfach, dass $\mathbb R$ genau die Menge aller Dezimalzahlen ist.

Insbesondere ist $\mathbb Q$ eine Teilmenge von $\mathbb R$ --  wir haben sogar unendlich viele Elemente hinzugef\"ugt, um von $\mathbb Q$ nach $\mathbb R$ zu kommen. Wie k\"onnen wir uns aber vorstellen, was genau bei der Erweiterung von $\mathbb Q$ nach $\mathbb R$ \emph{wirklich passiert}?  Das ist nicht ganz so elementar zu beschreiben. Wie wir gesehen haben, ist $\mathbb Q$ ziemlich dicht gepackt. Trotzdem ist $\mathbb Q$ noch \"ubers\"aht mit unendlich vielen, \emph{unendlich kleinen} ``L\"ochern'': $\mathbb Q$ ist unvollst\"andig. Zum Beispiel ist bei $\sqrt 2$ ein ``Loch'' in $\mathbb Q$, in folgendem Sinne: Wir k\"onnen $\sqrt 2= 1,41421...$ beliebig mit rationalen Zahlen ann\"ahern, durch die Zahlenfolge
$$1\;;\qquad 1,4\;; \qquad 1,41\;;\qquad 1,414\;;\qquad 1,4142\;;\qquad \text{usw.}$$
Da aber $\sqrt{2}$ selbst keine rationale Zahl ist, springen wir auf diese Weise in unendlich vielen Schritten ``aus $\mathbb Q$ heraus''. Und an der Stelle, wo wir aus $\mathbb Q$ herausfallen, ist gerade das Loch $\sqrt 2$.

Erst durch das Zusammenfassen aller Dezimalzahlen schlie\ss{}en wir alle L\"ocher, denn aus $\mathbb R$ k\"onnen wir nicht mehr auf die beschriebene Art und Weise herausspringen. Da $\mathbb R$ also keine L\"ocher hat, vollst\"andig ist und man sch\"on glatt darauf entlanggleiten kann, nennt man die reellen Zahlen auch das \emph{Kontinuum}.

\subsection*{Die \"Uberabz\"ahlbarkeit der reellen Zahlen -- \\ Cantors zweites Diagonalargument}
Wie wir gesehen haben geschieht etwas im \emph{unendlich Kleinen} beim \"Ubergang von $\mathbb Q$ nach $\mathbb R$. Dies hat auch gewaltige Auswirkungen auf das unendlich Gro\ss{}e. Denn mit einem weiteren Geniestreich konnte Cantor zeigen, dass die reellen Zahlen nicht abz\"ahlbar sind, dass also $\mathbb R$ m\"achtiger als $\mathbb N$ ist. Bevor wir uns sein wundersch\"ones Argument, sein \emph{zweites Diagonalargument}, genauer anschauen, erinnern wir uns an unsere Leitfrage vom Anfang, ob es irgendeine Struktur innerhalb der Aktualunendlichkeit gibt. Diese k\"onnen wir nun mit "Ja!" beantworten:
\begin{quote}
\emph{Es gibt verschiedene unendliche M\"achtigkeiten.}
\end{quote}
Die Bedeutung von Cantors Entdeckung l\"asst sich gar nicht genug w\"urdigen. Versuchen wir es mit den Worten von Bertrand Russell (1902):
\begin{quote}
``Die L\"osung der Schwierigkeiten, die fr\"uher das mathe\-matische Unendliche umgaben, ist wahrscheinlich die gr\"o\ss{}te Leistung, deren sich unser Zeitalter r\"uhmen kann."
\end{quote}

Wie hat Cantor gezeigt, dass $\mathbb R$ m\"achtiger als $\mathbb N$ ist? Dazu sind ja zwei Dinge nachzupr\"ufen: erstens dass $\mathbb N$ gleichm\"achtig zu einer Teilmenge von $\mathbb R$ ist, und zweitens dass $\mathbb N$ und $\mathbb R$ nicht gleichm\"achtig sind. Der erste Punkt ist leicht einzusehen: $\mathbb N$ ist n\"amlich selbst eine Teilmenge von $\mathbb R$ (wir erinnern uns dazu an $2 = 2,00000...$ und so weiter).

Es kommt also darauf an zu zeigen, dass $\mathbb N$ und $\mathbb R$ nicht gleichm\"achtig sind. Wie in unserem fr\"uheren Beispiel mit den endlichen Mengen $A = \{1,2,3\}$ und $B = \{a,b,c,d\}$ bleibt uns nichts anderes \"ubrig, als von allen Paarbildungen zu zeigen, dass immer mindestens ein Element von $\mathbb R$ \"ubrig bleibt.

Wie soll das m\"oglich sein? Im Beispiel waren 24 Paarbildungen zu pr\"ufen, aber nun gibt es doch unendlich viele davon! Dieses Problem k\"onnen wir durch folgende Strategie umgehen: Wir nehmen uns eine \emph{beliebige} Paarbildung zwischen $\mathbb N$ und $\mathbb R$ her und geben ein Element aus $\mathbb R$ direkt an, das keinen Partner aus $\mathbb N$ abbekommen hat. Wenn uns das gelingt, k\"onnen wir ganz sicher sein, dass es \emph{keine} Paarbildung ohne Rest gibt.

Cantors Idee, bei einer beliebigen Paarbildung zwischen $\mathbb N$ und $\mathbb R$ eine Dezimalzahl zu finden, die keinen Partner aus $\mathbb N$ hat, ist die folgende. Wir nehmen uns irgendeine Paarbildung zwischen $\mathbb N$ und $\mathbb R$ her, zum Beispiel 
$$
\begin{array}{ccc}
1 & \leftrightarrow & 0,3333...\\
2 & \leftrightarrow & 0,5432...\\
3 & \leftrightarrow & 0,6775...\\
4 & \leftrightarrow & 0,1010...\\
\vdots & \vdots & \vdots
\end{array}
$$
Wir konstruieren nun die gesuchte Dezimalzahl, indem wir die Diagonale auf der rechten Seite des Schaubilds entlanggehen und immer $1$ hinzuaddieren: vom Partner der $1$ nehmen wir uns die \emph{erste} Nachkommastelle und erh\"ohen sie um 1, dies wird die \emph{erste} Nachkommastelle unserer gesuchten Zahl, also $3+1=4$. Vom Partner der $2$ erh\"ohen wir die \emph{zweite} Nachkommastelle um $1$, was die \emph{zweite} Nachkommastelle unserer reellen Zahl ohne Partner wird, also $4+1=5$. Die \emph{dritte} Nachkommastelle wird die \emph{dritte} Nachkommastelle plus $1$ des Partners der $3$, die \emph{vierte} Nachkommastelle wird die \emph{vierte} Nachkommastelle plus $1$ des Partners der $4$, und so weiter. (Wenn wir dabei auf eine $9$ treffen, machen wir eine $0$ daraus.) Dieses Verfahren erzeugt eine Dezimalzahl, von der wir gleich sehen werden, dass sie keinen Partner haben kann. Machen wir uns das Verfahren am Beispiel klar:
$$
\begin{array}{ccc}
1 & \leftrightarrow & 0,\,\framebox[1.05\width]{\textbf{3}}\,333...\\
2 & \leftrightarrow & 0,5\,\framebox[1.05\width]{\textbf{4}}\,32...\\
3 & \leftrightarrow & 0,67\,\framebox[1.05\width]{\textbf{7}}\,5...\\
4 & \leftrightarrow & 0,101\,\framebox[1.05\width]{\textbf{0}}\,...\\
\vdots & \vdots& \vdots\\
? & \leftrightarrow & 0, \textbf{4581...}\\
\end{array}
$$
Die aus Cantors Diagonalverfahren resultierende Dezimalzahl ist hier also $0,4581...$\,. Kann nun die $1$ der Partner dieser Zahl sein? Nein, denn $0,4581...$ stimmt mit $0,3333...$ in der \emph{ersten} Nachkommastelle nicht \"uberein. Kann die $2$ der Partner sein? Nein, denn $0,4581...$ stimmt mit $0,5432...$ in der \emph{zweiten} Nachkommastelle nicht \"uberein. Und so weiter: die Zahl $0,4581...$ ist ja gerade so gemacht, dass sie mit jeder Dezimalzahl in der Liste auf mindestens eine Nachkommastelle nicht \"ubereinstimmt. 

Wir folgern, dass bei jeder beliebigen Paarbildung von $\mathbb N$ und $\mathbb R$ mindestens ein Element aus $\mathbb R$ keinen Partner hat. Eine Paarbildung ohne Rest ist also unm\"oglich. Damit ist $\mathbb R$ m\"achtiger als $\mathbb N$, oder auch: \emph{Die Menge der reellen Zahlen ist \"uberabz\"ahlbar.}

\"Ubrigens zeigt Cantors Argument sogar noch ein bisschen mehr: schon die Menge aller reellen Zahlen zwischen $0$ und $1$ ist \"uberabz\"ahlbar, denn oben haben wir ja nur Dezimalzahlen mit einer Null vor dem Komma betrachtet. Dies zeigt nochmal auf andere Weise, dass beim \"Ubergang von $\mathbb Q$ nach $\mathbb R$ das Entscheidende im unendlich Kleinen passiert.

\subsection*{Unendlich viele unendliche M\"achtig\-kei\-ten}
Wir wissen nun, dass es mindestens zwei verschiedene unendliche M\"achtig\-kei\-ten gibt, n\"amlich die von $\mathbb N$ und die von $\mathbb R$. Ohne Details zu nennen sei an dieser Stelle erw\"ahnt, dass es sogar eine ganze Abfolge unendlich vieler verschiedener M\"achtigkeiten gibt, die von Schritt zu Schritt immer m\"achtiger (also "gr\"o\ss{}er") werden. Damit k\"onnen wir festhalten:
\begin{quote}
\emph{Es gibt unendlich viele verschiedene unendliche M\"achtigkeiten.}
\end{quote}

\subsection*{Die Feinstruktur der Aktualunendlichkeit -- \\ Cantors Kontinuumshypothese}
In diesem letzten Abschnitt wollen wir noch zwei Resultate zur genaueren Struktur der Aktualunendlichkeit ansprechen, wiederum ohne zu sehr ins Detail zu gehen. Dabei erreichen wir prinzipielle Grenzen der Erkenntnis!

Zun\"achst dr\"angt sich uns die Frage auf, \emph{wie viele} verschiedene unendliche M\"achtigkeiten es denn gibt. Gibt es abz\"ahlbar viele, so viele wie es reelle Zahlen gibt, oder sogar noch mehr? 

Die Antwort darauf ist verbl\"uffend: Die Gesamtheit aller verschiedenen unendlichen M\"achtigkeiten ist so extrem gro\ss{}, dass man sie gar nicht alle gleich\-zeitig zu einer Menge zusammenfassen kann, ohne  verheerende Widerspr\"uche zu erzeugen. Schon bei der Zusammenfassung aller nat\"urlichen Zahlen zu einer Menge haben wir gesehen, dass zumindest Vorsicht angebracht ist -- wobei heute aber angenommen wird, dass dies unproblematisch ist. Das Zusammenfassen verschiedener Dinge zu einem Ganzen hat aber seine Grenzen -- und diese ist bei der Gesamtheit aller unendlichen M\"achtigkeiten erreicht.

Noch viel interessanter wird es aber in der N\"ahe der M\"achtigkeit der nat\"urlichen Zahlen $\mathbb N$. Wir wissen schon, dass dies die kleinste unendliche M\"achtigkeit ist. Au\ss{}erdem wissen wir, dass die Menge der reellen Zahlen $\mathbb R$, das Kontinuum, m\"achtiger als $\mathbb N$ ist. Deshalb ist es ganz nat\"urlich zu fragen:
\begin{quote}
\emph{Gibt es eine M\"achtigkeit \textbf{zwischen} denen von $\mathbb N$ und $\mathbb R$\,?}
\end{quote}
Anders gefragt: Ist jede unendliche Teilmenge von $\mathbb R$ entweder so m\"achtig wie $\mathbb N$ oder gleichm\"achtig zu $\mathbb R$\,?

Cantor hat sich diese fundamentale Frage als Erster gestellt. Jedoch konnte er sie weder positiv noch negativ beantworten. Er war aber fest davon \"uberzeugt, dass es \emph{keine} solche M\"achtigkeit gibt. 1878 ver\"offentlichte er diese Vermutung, die als \emph{Kontinuumshypothese} in die Geschichte eingegangen ist.

F\"ur viele Jahrzehnte war die Antwort auf die Frage unbekannt, Mathematiker und Philosophen konnten die Kontinuumshypothese weder beweisen noch widerlegen. David Hilbert stellte im Jahr 1900 eine Liste mit den zu der Zeit 23 wichtigsten ungel\"osten Problemen der Mathematik auf (ver\-gleich\-bar mit den sieben \emph{Milleniumsproblemen} aus dem Jahr 2000), bei denen die Kontinuumshypothese sogar an erster Stelle stand. 

Die letztendliche Antwort, die von Kurt G\"odel (1938) und Paul Cohen (1960) gefunden wurde, ist sensationell.

Schauen wir uns die Frage nochmal etwas genauer an. Es geht darum, ob Cantors Hypothese \emph{wahr oder falsch} ist. Dabei setzen wir also voraus, dass die Hypothese wirklich \emph{entweder wahr oder falsch} ist. Was sollte sie auch sonst sein? Die Annahme, \emph{dass} jede Aussage entweder wahr oder falsch ist, ist so etwas wie das Grundwerkzeug f\"ur jeden Mathematiker, vergleichbar mit (frei nach Hilbert) dem Fernrohr f\"ur den Astronomen oder dem Reagenzglas f\"ur den Chemiker. Es gilt aber folgendes:

\begin{quote}
\emph{Die Kontinuumshypothese ist weder beweisbar noch widerlegbar, weder wahr noch falsch, sie ist \emph{unentscheidbar}.}
\end{quote}
Das war nun wirklich nicht was man erwartet h\"atte. Wie kann das sein, wie kann etwas weder wahr noch falsch sein? 

Dazu muss man sich klar machen, was es \"uberhaupt hei\ss{}t, dass eine mathematische Aussage beweisbar bzw. widerlegbar ist. Eine Aussage zu \emph{beweisen} bedeutet n\"amlich, sie logisch aus Annahmen abzuleiten, die nicht mehr begr\"undet werden -- den Axiomen. Andererseits bedeutet eine Aussage zu widerlegen, deren Gegenteil aus den Axiomen logisch abzuleiten. 

Ob also eine Aussage \emph{entscheidbar} (das hei\ss{}t entweder beweisbar oder widerlegbar) ist, h\"angt von zwei Ding\-en ab: den Axiomen und den lo\-gi\-schen Ableitungsregeln. Es ist aber ganz und gar nicht klar, dass jede Aussage innerhalb eines Axiomensystems auch wirklich entscheid\-bar ist. Die Mathematiker haben dies bis 1931 einfach stillschweigend angenommen.

In unserem Fall sind die zugrundeliegenden Axiome gerade die \emph{Axiome der Mengenlehre}. Ein Axiom der Mengenlehre haben wir schon kennen\-ge\-lernt, n\"amlich das \emph{Unendlichkeitsaxiom}, welches besagt, dass es eine Menge gibt, die alle nat\"ur\-lichen Zahlen enth\"alt. Weitere Axiome der Mengenlehre sind zum Beispiel, dass es eine Menge gibt, die keine Elemente enth\"alt (die \emph{leere Menge}) und dass man zwei Mengen immer zu einer neuen Menge vereinigen kann. Dass es eine Menge gibt, die alle reellen Zahlen enth\"alt, dass man also alle reellen Zahlen wirklich zu einer Menge zusammenfassen darf, ist eine Konsequenz aus dem \emph{Potenzmengenaxiom}. Wir haben auch schon angesprochen, was die Axiome der Mengenlehre unter anderem nicht erlauben: n\"amlich alle verschiedenen unendlichen M\"achtigkeiten zu einer Menge zusammenzufassen.

Die Unentscheidbarkeit der Kontinuumshypothese ist nun folgender\-ma\ss{}\-en zu verstehen: Man kann \emph{beweisen} (G\"odel), dass man das Gegenteil der Hypothese nicht aus den Axiomen der Mengenlehre ableiten kann. Andererseits kann man aber \emph{auch beweisen} (Cohen), dass man die Hypothese selbst auch nicht aus den Axiomen der Mengenlehre ableiten kann. 

Das bedeutet, dass die Kontinuumshypothese \emph{unabh\"angig} von den Axiomen der Mengenlehre ist. Man hat nun zumindest prinzipiell die Wahl, ob man sie als \emph{neues Axiom} zu den anderen Axiomen der Mengenlehre hinzunimmt, ob man ihr Gegenteil als neues Axiom mit hinzunimmt -- oder ob man sie einfach ignoriert.

Die Unentscheidbarkeit von Cantors Hypothese war das erste relevante Beispiel f\"ur ein grunds\"atzliches Problem der Mathematik: G\"odel zeigte n\"amlich 1931 in seinen \emph{Unvollst\"andigkeitss\"atzen} unter anderem, dass die allermeisten f\"ur die Mathematik relevanten Axiomensysteme un\-ent\-scheid\-bare Aussagen enthalten. Genauer betrifft dies alle Axiomensysteme, deren Informationsgehalt ausreicht, um das \"ubliche Addieren und Multiplizieren mit nat\"urlichen Zahlen zu beschreiben. Das sind wahrlich keine gro\ss{}en Anforderungen. G\"odels Resultate haben die Mathematik in ihren Grundfesten ersch\"uttert.

Bis zum Nachweis der Unentscheidbarkeit der Kontinuumshypothese im Jahr 1960 kannte man aber nur ''k\"unstliche" unentscheidbare Aussagen. (Ein Beispiel ist die Aussage "Diese Aussage ist falsch". Nehmen wir n\"amlich an, die Aussage "Diese Aussage ist falsch" w\"are wahr, dann w\"are das im Widerspruch zur Aussage selbst, die ja aussagt, dass die Aussage falsch ist. Nehmen wir aber umgekehrt an, die Aussage "Diese Aussage ist falsch" w\"are falsch, dann ergibt das auch einen Widerspruch, da das Gegenteil der Aussage n\"amlich gerade hei\ss{}t, die Aussage "Diese Aussage ist falsch" w\"are wahr. Somit ist diese Aussage weder beweisbar noch widerlegbar. Alles klar?)

Sowohl das heute akzeptierte Axiomensystem der Mengenlehre, als auch die g\"ultigen logischen Ableitungsregeln sind von Menschenhand gemacht und deshalb keineswegs bis in alle Ewigkeit in Stein gemei\ss{}elt. Bis heute ist eine m\"ogliche \"Anderung des klassischen Axiomensystems der Mengenlehre -- hin zur Entscheidbarkeit der Mengenlehre -- Gegenstand der mathe\-ma\-tischen Forschung. Es ist noch nicht abzusehen, in welche Richtung sich diese Geschichte entwickeln wird.

\subsection*{Weiterf\"uhrende Literatur}
\begin{itemize}
\item Spektrum der Wissenschaft Highlights 2/2013: \emph{Das Unendliche.}
\item H. Heuser. \emph{Unendlichkeiten: Nachrichten aus dem Grand Canyon des Geistes}. Teubner (2007).
\item D. R. Hofstadter. \emph{G\"odel, Escher, Bach: Ein endloses geflochtenes Band}. dtv (1992).
\item O. Deiser. \emph{Einf\"uhrung in die Mengenlehre: Die Mengenlehre Georg Cantors und ihre Axiomatisierung durch Ernst Zermelo.} Springer-Verlag (2010).
\end{itemize}

\end{document}